\documentclass[a4paper]{article}
\usepackage{amsmath}
\usepackage{bbm}        % nicer symbols for sets
\usepackage{mma}
\usepackage{subfigure}
\usepackage[pdftex]{graphicx}

\def\D#1{D_{\!#1}}
\def\S#1{S_{\!#1}}
\def\N{\mathbbm N}

\def\Q{\mathbbm Q}
\def\R{\mathbbm R}

\def\d{\partial}
\def\rmd{\,\mathrm{d}}

\def\blank{\underline{\;\;}}                                  % Blank[]
\def\pow{\mathbin{\raisebox{-2.5pt}{\hbox{\large$\hat{}$}}}}  % Power x^y
\def\curl{\operatorname{curl}}                                % curl 
\def\sumoverallelements{\sum_{T \in \mathcal{T}_h}}
\def\sumoverallfacets{\sum_{F \in \mathcal{F}_h}}
\def\normal{\nu}
\def\averageleft{\{\!\!\{}
\def\averageright{\}\!\!\}}
\def\jumpleft{[\![}
\def\jumpright{]\!]}

\begin{document}

\title{Computer Algebra meets Finite Elements: an Efficient Implementation for Maxwell's Equations\footnote{This article is part of the volume U.~Langer and P.~Paule (eds.) \emph{Numerical and Symbolic Scientific Computing: Progress and Prospects} in the series \emph{Texts \& Monographs in Symbolic Computation}, ISBN 978-3-7091-0793-5. The original publication is available at www.springerlink.com, DOI 10.1007/978-3-7091-0794-2\_6.}}
% Efficient Implementation of an Higher Order DG Method for the Time-Domain Maxwell's Equations}

\author{
     Christoph~Koutschan\thanks{supported by the Austrian Science Fund (FWF): SFB F013 and P20162-N18,
       and partially by NFS-DMS 0070567 as a postdoctoral fellow.}\\
     Research Institute for Symbolic Computation\\
     Johannes Kepler University\\
     Linz, Austria
\and Christoph~Lehrenfeld\\
     Institut f\"ur Geometrie\\
     und Praktische Mathematik\\ 
     RWTH Aachen, Germany
\and Joachim~Sch\"oberl\\ 
     Center for Computational\\ 
     Engineering Science\\ 
     RWTH Aachen, Germany
}

\maketitle

\begin{abstract}
We consider the numerical discretization of the time-domain Maxwell's
equations with an energy-conserving discontinuous Galerkin finite
element formulation. This particular formulation allows for higher
order approximations of the electric and magnetic field. Special
emphasis is placed on an efficient implementation which is achieved by
taking advantage of recurrence properties and the tensor-product
structure of the chosen shape functions. These recurrences have been
derived symbolically with computer algebra methods reminiscent of the
holonomic systems approach.
\end{abstract}

\section{Introduction}
This paper is dedicated to a successful cooperation between symbolic
computation and numerical analysis.  The goal is to simulate the
propagation of electromagnetic waves using finite element methods
(FEM).  Such simulations play an important role for constructing
antennas, electric circuit boards, bodyworks, and many other devices
where electromagnetic radiation is involved. The numerical simulation
of such physical phenomena helps to optimize the shape of components
and saves the engineer from doing a long and expensive series of
experiments.

Finite element methods serve to approximate the solution of partial
differential equations on a given domain~$\Omega\subseteq\R^d$ subject
to certain constraints (e.g., boundary conditions).  The
domain~$\Omega$ is partitioned into small elements (typically
triangles or tetrahedra) and the solution is approximated on each
element by means of certain shape functions. In our application we
deal with Maxwell's equations which relate the magnetic and the
electric field. In Section~\ref{secFormulation} we describe how the
problem can be discretized using FEM and in
Section~\ref{secImplementation} we give the details concerning an
efficient implementation.

An important ingredient for the fast execution of some operations in
the FEM are certain difference-differential relations that were
derived with computer algebra methods.  The methods that we employ,
originate in Zeilberger's holonomic systems approach
\cite{Zeilberger90,Chyzak00,Koutschan09} whose basic idea is to define
functions and sequences in terms of differential equations and
recurrence equations plus initial values (these equations have to be
linear with polynomial coefficients).  Luckily the shape functions
used in the chosen FEM discretization fit into the holonomic framework since
they are defined in terms of orthogonal polynomials.
Section~\ref{secSymbolics} explains how the desired relations have
been computed.

\section{FEM formulation of Maxwell's equations}\label{secFormulation}
In order to describe electromagnetic wave propagation problems,
we consider the loss-free time-domain Maxwell's equations 
\begin{eqnarray*}
  \varepsilon\frac{\d E}{\d t}  & = & \curl H,\\
  \mu\frac{\d H}{\d t} & = & -\curl E,
\end{eqnarray*}
subject to appropriate initial and boundary conditions.
Here $E=E(x,t)$ denotes the electric and $H=H(x,t)$ the magnetic field strength
(with $x=(x_1,x_2,x_3)$ the space variables and $t$ the time),
and~$\varepsilon$ and $\mu>0$ are the permittivity and 
the permeability, respectively. When discretizing these equations with the 
finite element method, we go over to a weak formulation by multiplying 
both equations with test functions~$e(x)$ and~$h(x)$ and integrating over 
the whole domain~$\Omega\subset\R^3$. The solution of the Maxwell's equations then has to
fulfill the conditions
\begin{equation}\label{weakform}
  \begin{array}{rcl}
    \displaystyle\frac{\d}{\d t}(\varepsilon E,e)_\Omega & = & (\curl H,e)_\Omega,\\[2ex]
    \displaystyle\frac{\d}{\d t}(\mu H,h)_\Omega & = & -(\curl E,h)_\Omega
  \end{array}
\end{equation}
for all test functions~$e$ and~$h$,
where $(\cdot,\cdot)_\Omega$ is the short notation for the $L^2(\Omega)$ inner product 
$(a,b)_\Omega=\int_\Omega ab\rmd x$.
Then we replace both the magnetic and electric field as well as the
test functions by finite-dimensional approximations 
on a triangulation~$\mathcal{T}_h$ of the domain~$\Omega$. 
Herein $h$ denotes some characteristic length of the elements in~$\mathcal{T}_h$
(not to be confused with the test function~$h$).

Conforming finite elements ensure that the finite-dimensional 
approximations are within a space which is appropriate for the partial 
differential equations under consideration. For Maxwell's equations this 
space is $H(\curl,\Omega)$ which demands tangential components 
to be continuous across element interfaces. 
The discontinuous Galerkin finite element method (DG) neglects this 
conformity condition when building up a discrete basis for the 
approximation, but instead has to incorporate stabilization terms to achieve a
consistent and stable formulation. This is normally done by applying 
integration by parts and replacing fluxes at element 
boundaries with \emph{numerical fluxes}
\cite{ArnoldBrezziCockburnMarini, psm02, hps04, HesthavenWarburton03}.
The latter approach has the major advantage that the mass matrices $M_{\varepsilon}$ and $M_{\mu}$, i.e.,
the matrices that arise when discretizing $(\varepsilon E,e)_{\Omega}$ and 
$(\mu H,h)_{\Omega}$, respectively, are block-diagonal which makes the 
application of their inverses computationally more efficient. 

We consider the approximation space 
\[
  V_h^k = \left\{ v \in\left(L^2(\Omega)\right)^3: v|_T \in \left(\mathcal{P}^{k}(T)\right)^3\ \forall T \in \mathcal{T}_h \right\}
\]
that consists of functions which are piecewise polynomial up to degree~$k$. 
By integration by parts of \eqref{weakform} on each element~$T\in\mathcal{T}_h$, and by adding a 
consistent stabilization term on all element boundaries we get
(again for all test functions~$e$ and~$h$)
\begin{eqnarray*}
  \frac{\d}{\d t} \sumoverallelements (\varepsilon E,e)_{T} &=& 
    \sumoverallelements \big( (H, \curl e)_T + (H^* \times \normal,e)_{\d T} \big),\\
  \frac{\d }{\d t} \sumoverallelements (\mu H,h)_{T} &=& 
    \sumoverallelements \big( -(\curl {E},h)_T + (E^*-E , h \times \normal)_{\d T} \big), 
\end{eqnarray*}
where $\normal$ denotes the outer normal on each element boundary and $H^*$, 
$E^*$ are the numerical fluxes.
The properties of different DG formulations mainly depend on 
the choice of the numerical fluxes. As all derivatives 
are now shifted to the electric field~$E$ and the according test 
functions~$e$, it is reasonable to approximate the electric 
field of one degree higher than the magnetic field. 
So we choose the approximation 
spaces $V_h^{k+1}$ for~$E$ and~$e$ and $V_h^{k}$ for~$H$ and~$h$.

\subsection{Numerical flux}
Several choices for the numerical flux are used in practice. 
Our goal here is to derive a numerical flux which ensures that 
the numerical approximation fulfills the following two important 
properties which are already fulfilled on the continuous level:
\begin{enumerate}
 \item conservation of the energy 
	$\frac12(\varepsilon E,E)_{\Omega} + \frac12(\mu H,H)_{\Omega}$
 \item non-existence of \emph{spurious modes}
\end{enumerate}
On the one hand using dissipative fluxes avoids spurious modes
and is often used, but as it introduces dissipation, the energy 
of the system is not 
conserved. On the other hand the standard approach for energy 
conserving methods is the so called \emph{central flux}. Its mayor
disadvantage is, that it introduces non-physical modes, 
spurious modes. 

Nevertheless
we start with this approach to derive the \emph{stabilized central flux}
formulation which gets rid of both problems. 
A more extensive discussion of numerical fluxes (including
the stabilized central flux) for Maxwell's equations can be found in 
\cite[8.2]{HesthavenWarburton07}. 

The central flux takes the averaged values of neighboring elements
for the numerical flux, i.e., $H^* = \averageleft H \averageright$ and
$E^* = \averageleft E \averageright$ with
$\averageleft\cdot\averageright$ denoting the averaging operator, and
ends up with a semi-discrete system of the form
\begin{equation}\label{semidiscrete}
  \frac{\d}{\d t}\begin{pmatrix}M_{\varepsilon} & \\ & M_{\mu} \end{pmatrix}\begin{pmatrix}E \\ H\end{pmatrix}=
  \begin{pmatrix} & -C_h^T\\ C_h& \end{pmatrix}\begin{pmatrix}E \\ H\end{pmatrix}
\end{equation}
where $C_h$ denotes the discrete $\curl$ operator stemming from the
central flux formulation. The matrix on the left side is symmetric and positive 
definite whereas the matrix on the right side is antisymmetric. Then the evolution matrix for 
the modified unknowns $(M_{\varepsilon}^{\frac12}E, M_{\mu}^{\frac12}H)^T$ is 
also antisymmetric and thus the proposed energy is conserved.
Nevertheless this matrix has a lot of eigenvalues close to zero which 
correspond to the discretization, but not to the physical behavior of the 
system.
To motivate the modification which will stabilize the formulation, let 
us have a brief look at the problem in frequency domain, i.e., for time-harmonic 
electric and magnetic fields. Then
the discrete problem in frequency domain reads (with frequency~$\omega$):
\begin{equation}\label{discretefrequencydomain}
  0 = (i \omega)^2 (M_\varepsilon E,e) + (M_\mu^{-1} C_h E, C_h e).
\end{equation}
The problem with non-physical zero eigenvalues now manifests
in $(C_h E, C_h e)$ being only positive semidefinite. We overcome this 
issue by adding a stabilization bilinear form $S(E,e)$ to 
\eqref{discretefrequencydomain}  as proposed in \cite{HesthavenWarburton07}.
\[
  S(E,e) := \sumoverallfacets \frac{\alpha}{h} (\jumpleft E \jumpright \times \normal, \jumpleft e \jumpright \times \normal)_F
\]
with $\alpha>0$, where $\mathcal{F}_h$ is the union of all element boundaries and 
$\jumpleft\cdot\jumpright$ denotes the jump operator, i.e., the difference 
between values of adjacent elements. This stabilization bilinearform eliminates 
the nontrivial kernel of $C_h$ and is consistent as 
$\jumpleft E \jumpright \times \normal$ is zero for the exact solution.
Before we can translate the formulation back to the time domain, 
we introduce a new variable which is defined as
\[
  H^F := \frac{(\jumpleft E \jumpright \times \normal) \alpha}{i \omega h}
\]
The new unknown~$H^F$ is also piecewise polynomial on each face. 

If we go back to the time-domain formulation we end up with the following
formulation (note that relations between $\jumpleft \cdot \jumpright$ and
$\averageleft \cdot \averageright$ were used):
\begin{eqnarray*}
  \frac{\d}{\d t} \sumoverallelements (\varepsilon E,e)_{T} & = &
    \sumoverallelements \big( (H, \curl e)_T  + (\averageleft H \averageright \times \normal,e)_{\d T} \big)+\\
  &&  \sumoverallfacets  (H^F \times \normal,\jumpleft e \jumpright)_F, \\ 
  \frac{\d}{\d t} \sumoverallelements (\mu H,h)_{T} & = & 
    \sumoverallelements \big( -(h, \curl E)_T  + (\textstyle\frac12 \jumpleft E \jumpright \times \normal,h)_{\d T}\big),\\ 
  \frac{\d}{\d t} \sumoverallfacets \frac{\alpha}{h} (H^F,h^F)_F & = &
    \sumoverallfacets (\jumpleft E \jumpright \times \normal,h^F)_F. 
\end{eqnarray*}
For $p$-robust behavior $\alpha$ should scale with $p^2$, where $p$ 
is the polynomial degree. This is motivated by the symmetric interior 
penalty method for elliptic equations (see e.g. \cite{ArnoldBrezziCockburnMarini}) 
where a scaling of $\alpha$ with $p^2$ in the bilinearform $S$ is 
necessary for stability to dominate over some terms stemming 
from inverse inequalities which scale with $p^2$ 
(see also \cite{HesthavenWarburton03}).

We again achieve a system of the form \eqref{semidiscrete} where the 
vector $H$ now consists of element and face unknowns and the matrix
representing the discrete $\curl$ operator is the stabilized 
central flux $\curl$ operator now. 
Thus we conclude that the method now conserves energy, and 
spurious modes, introduced by the central flux, vanish.

\subsection{Numerical Examples (Spherical Vacuum Resonator)}
\begin{figure}[ht]
  \begin{center}
  \includegraphics[width=12cm]{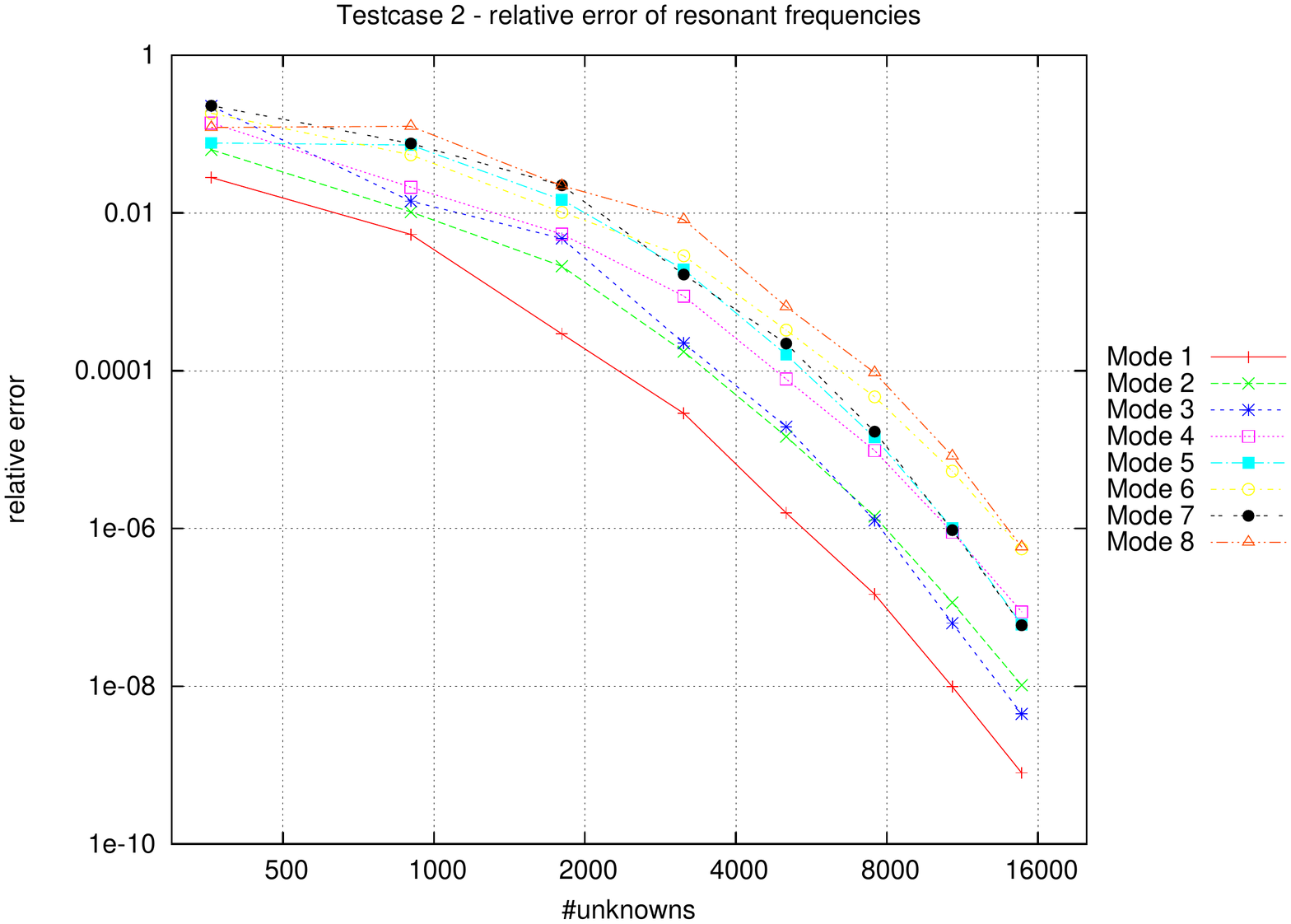}
  \end{center}
  \vspace{-10mm}
  \caption{Convergence of the resonance frequencies after $p$-refinement}\label{fig:resonator}
\end{figure}

We consider a spherical domain $\Omega := \{ x \in \R^3: \Vert x \Vert_2 \leq 1\}$
and the frequency domain formulation of the Maxwell's equations subject to perfect electrical boundary conditions 
\begin{eqnarray*}
\begin{array}{ccc}
\left. \begin{array}{rcl}
  i \omega \varepsilon E  & = & \curl H,\\
  i \omega \mu H & = & -\curl E,
\end{array} \right\} 
\quad & \text{on} & \quad \Omega, \\
E \times \normal = 0 & \text{on} & \quad \partial \Omega,
% the boundary)
\end{array}
\end{eqnarray*}
% \begin{equation}
% \end{equation}
To demonstrate the opportunities of higher order discretizations we
consider a coarse mesh consisting of 30~elements and increase the
polynomial degree to increase the spatial resolution. We are
interested in the error of the eight smallest resonance
frequencies. Therefore we compare the eigenvalues of the numerical
discretization with those of a reference solution.  In
Figure~\ref{fig:resonator} we observe the expected exponential
convergence of the method.

\section{Computational aspects}\label{secImplementation}
As the spatial discretization conserves energy, we consider symplectic
time integration methods which conserve the energy on a time-discrete
level. The simplest one is the \emph{symplectic Euler} method which 
discretizes the semi-discrete system \eqref{semidiscrete} in the following way:
\begin{eqnarray*}
 H^{n+1} & = & H^n + \Delta t \ M_{\mu}^{-1} C_h E^n \\
 E^{n+1} & = & E^n - \Delta t \ M_{\varepsilon}^{-1} C_h^T H^{n+1}
\end{eqnarray*}
with the stability condition 
\[
\Delta t \leq 2  \left( \rho( M_\mu^{-\frac12} C_h M_e^{-1} C_h^T M_\mu^{-\frac12} ) \right)^{-1}
\]
The matrix $M_\mu^{-\frac12} C_h M_e^{-1} C_h^T M_\mu^{-\frac12}$ is
symmetric and the spectral radius $\rho$ can be estimated once by an
iterative method like the power iteration. \\ When shifting the
electric or the magnetic field by a half time-step we can reconstruct
the well-known \emph{leap frog} method. Nevertheless for our
considerations it is less important which time integration scheme is
used as long as it is explicit. The matrix multiplications with $C_h$
and $C_h^T$ (see Section~\ref{secCurls}) as well as with
$M_{\mu}^{-1}$ and $M_{\varepsilon}^{-1}$ (see
Section~\ref{secMassMatrices}) decide about the computational
efficiency of an implementation.

The advantage of discontinuous Galerkin methods becomes evident now. 
The mass matrices can be inverted in an element by element fashion and also
the discrete $\curl$ operations only need information of (element-)local and 
adjacent degrees of freedom, which allows for straightforward parallelization.
Element matrices such as mass matrices and the discrete $\curl$ operation 
can be stored once and applied at each time step. This is how far one comes 
just because of the formulation itself. 

With appropriate choices for the local shape functions we can use advanced 
techniques to execute those operations with a lower complexity than local 
matrix-vector multiplications. Furthermore we don't even have to store the element
matrices, s.t. the techniques presented below are also much more memory-efficient.

The following ingredients are essential for the techniques proposed below, 
which enhance the implementation of the DG method:
\begin{enumerate}
 \item Definition of an $L^2$-orthogonal basis of polynomial shape functions in tensor-product
form\footnote{these are polynomials which are products of univariate
polynomials}  on a reference element $\hat{T}$
 \item Use of $\curl$-conforming (\emph{covariant}) transformation for evaluations on the physical element~$T$
 \item Use of recurrences for the polynomial shape functions to
evaluate gradients and $\curl$s
 \item Use of tensor-product structure to evaluate
traces\footnote{values at a boundary}
\end{enumerate}

\subsection{Local shape functions}
For stability and fast computability we choose the $L_2$-orthogonal
Dubiner basis \cite{dubiner91, Karniadakis}. Here, the basis functions
on the reference element are constructed in a tensor-product form of
Jacobi polynomials $P_i^{(\alpha,\beta)}$ for each spatial component
(note that the Legendre polynomials~$P_i=P_i^{(0,0)}$ are just a
special case).  For example, on the reference triangle spanned by the
points $(0,0)$, $(1,0)$ and $(0,1)$ the shape functions take the form
\begin{equation}\label{phi2D}
  \varphi_{i,j}(x,y) = P_i\left(\textstyle\frac{2y}{1-x}-1\right) \cdot (1-x)^i \cdot P_j^{(2i+1,0)}(2x-1).
\end{equation}

They are orthogonal on the reference triangle, and gradients can be evaluated 
by means of recurrence relations as demonstrated in Section~\ref{secGradients}. 
Due to the tensor-product form traces can be evaluated very fast, see
Section~\ref{secTraces}.

\subsection{Discrete $\curl$ operations}\label{secCurls}
At each time step 
% for functions $H$ and $H^F$ which is given in modal representation, 
we have to evaluate terms like $(H, \curl e)_T$ on each element~$T$ and 
$(\averageleft H \averageright \times\normal,\jumpleft e \jumpright)_F$ on each face~$F$.
Similar expressions have to be evaluated for the electric field~$E$.

\subsubsection{Covariant transformation}
Let $\Phi: \hat{T} \rightarrow T$ be a diffeomorphic mapping from the reference element to some physical element~$T$. 
Then the covariant transformation of a function $\hat{u}$ defined on the
reference element $\hat{T}$ is
\begin{equation*}
u := (F^{-1})^T \hat{u} \circ \Phi^{-1} \qquad \mbox{ with } \qquad F = \nabla
\Phi.
\end{equation*}
% with $J = \det(F'_T)$, implies $\vect{u} \in H(\divergence,T)$. The divergence of $\hat{\vect{u}}$ is transformed even simpler:
% \begin{equation*}
%  \divergence(\vect{u}) := J^{-1} \ \divergence_{\hat{\vect{x}}}(\hat{\vect{u}} \circ F_T^{-1})
% \end{equation*}
If we define the shape functions on the mapped elements as the
covariant transformed shape functions on the reference element, 
then the tangential component on the mapped element depends only on the tangential
component of the reference element.
The transformation
is called $\curl$-conforming as it ensures that for any function $\hat{u} \in
H(\curl,\hat{\Omega})$ the covariant transformed function~$u$ lies in
$H(\curl,\Omega)$. Furthermore it preserves certain integrals, s.t.
the following relations hold for the covariant transformations $H,e \in
H(\curl,T)$ of $\hat{H},\hat{e} \in H(\curl,\hat{T})$:
\begin{eqnarray*}
\left|\int_T H\curl e \rmd x\right| &=& \left|\int_{\hat{T}} \hat{H} \curl \hat{e} \rmd x\right|, \\
\left|\int_{\d T} (H \times \normal) e \rmd s\right| &=& \left|\int_{\d \hat{T}} (\hat{H} \times \normal) \hat{e} \rmd s\right|. 
\end{eqnarray*}
This means that the integrals of these forms appearing in the formulation 
are independent of the  geometry of the particular elements. The matrices can be 
computed once on the reference element. This trick was published in \cite{cohen06}.

\subsubsection{Evaluating gradients}\label{secGradients}
For computing $\curl$s it is sufficient to evaluate gradients, since
the $\curl$ is a certain linear combination of derivatives.  We write
the corresponding function~$\hat{E}$ in modal representation, i.e.,
\[
  \hat{E}=\sum_\alpha a_\alpha\varphi_\alpha,\ a_\alpha\in\R^3,
\]
where the sum ranges over the finite collection of (scalar) shape
functions defined on the reference element (in 2D the multi-index
$\alpha$ is $(i,j)$ and in 3D $\alpha=(i,j,k)$). With the use of the
covariant transformation, we just have to consider the integral on the
reference element~$\hat{T}$:
\[
  \int_{\hat{T}} \hat{h} \curl \hat{E} \rmd x.
\]
The idea is now to take advantage of recurrence relations between 
derivatives of Jacobi polynomials and Jacobi polynomials itself. 
We aim for an operation which gives the coefficients $b_\alpha \in \R^3$ representing the gradient
\[
  \nabla \hat{E} = \sum_\alpha b_\alpha \varphi_\alpha.
\]
Then $L^2$-orthogonality can be used to evaluate the complete integral 
very fast.

For ease of presentation let's consider the far more easy case of
evaluating the derivative of a scalar one-dimensional function 
$v(x) = \sum_{i=0}^n v_i P_i(x),\ v_i\in\R$ 
given in a modal basis of Legendre polynomials $P_i$, which fulfill the relation
\begin{equation}\label{rewriteLegendre}
  P'_{i+1}(x) = P'_{i-1}(x) + (2i+1)P_i(x).
\end{equation}
Then the problem is to find the modal representation of
\[
  v'(x) = \sum_{i=0}^n v_i P_i^\prime(x) = \sum_{i=0}^{n-1} w_i P_i(x).
\]
Let's show the first step, i.e., how we get the highest order
coefficient $w_{n-1}$:
\begin{eqnarray*}
v^\prime(x) 
        & = & \sum_{i=0}^{n} v_i P_i^\prime(x) 
	= \sum_{i=0}^{n-1} v_i P_i^\prime(x) + v_n P_n^\prime(x) \\
        & = & \sum_{i=0}^{n-1} v_i P_i^\prime(x) + v_n P_{n-2}^\prime(x) + v_n (2n-1) P_{n-1} \\
        & = & \sum_{i=0}^{n-1} \tilde v_i P_i^\prime(x) + w_{n-1} P_{n-1}(x)
\end{eqnarray*}
where we used the recurrence relation \eqref{rewriteLegendre} for $P_n^\prime(x)$ and thus get $w_{n-1} = v_n (2n-1)$. 
For the remaining polynomial $\sum_{i=0}^{n-1} \tilde v_i P_i^\prime(x)$ of degree $n-1$ we can apply 
the same procedure to get $w_{n-2}$. This can be continued until also $w_0$ and thereby the complete polynomial
representation $\sum_{i=0}^{n-1} w_i P_i(x)$ of $v^\prime (x)$ is determined.

An efficient C++ implementation of this procedure was achieved by 
\emph{template meta-programming}, where the compiler can generate 
optimized code for all elements up to an a priori chosen maximal
polynomial order.

The same basically also works in three dimensions with Jacobi
polynomials, but the relations are far more complicated, see
Section~\ref{secSymbolics}, and need 3 nested loops.

The overall costs for the evaluation of the
element $\curl$ integral scales \emph{linearly} with the number
of unknowns $N$ on one element which is much better than the 
matrix-vector multiplication which already has complexity 
$\mathcal{O}(N^2)$.

\subsubsection{Evaluating traces}\label{secTraces}
The boundary integrals that have to be evaluated can make use
of the tensor-product form to evaluate traces. Again we don't want 
those traces to be evaluated pointwise but in a modal sense and 
recurrences for the Jacobi polynomials make the transformation from 
volume element shape functions to face shape functions with 
$\mathcal{O}(N)$ operations possible. The procedure therefore 
is similar to the evaluation of the gradient in the previous section.

\subsection{Mass matrix operations}\label{secMassMatrices}
So far we dealt only with the discrete $\curl$ operations. So the only 
thing that is left to talk about is the application of the inverse mass matrices.
Due to the covariant transformation we have
\begin{eqnarray}\label{massintegral}
 ((M_{\varepsilon})_{\alpha,\beta})_{l,m} & = &
 \int_T \varepsilon\,(\varphi_\alpha e_l^T)\,(\varphi_\beta e_m) \rmd x \nonumber\\
 & = & \int_{\hat{T}} |\det(F)|\,\varepsilon\,(\hat{\varphi}_\alpha e_l^T) F^{-1} (F^{-1})^T (\hat{\varphi}_\beta e_m) \rmd x
\end{eqnarray}
with $\varphi_\alpha$ denoting the scalar-valued shape functions and $e_n$ the $n$-th unit vector.
Note also the block structure of~$M_\varepsilon$ that is indicated by the above notation.
In some FEM applications, symbolic methods related to those described in Section~\ref{secSymbolics},
can be used to prove the sparseness of the corresponding system matrix, see~\cite{Pillwein08}.

\subsubsection{Flat elements}
Let's assume the material parameters~$\varepsilon$ and~$\mu$ are 
piecewise constant and the elements are flat, i.e., 
$\nabla\Phi = F = const$ on each element. Then the integral~\eqref{massintegral} simplifies to
\[
 \int_T \varepsilon\,(\varphi_\alpha e_l^T)\,(\varphi_\beta e_m) \rmd x= 
 \vert \det(F) \vert\, \varepsilon\, (F^{-1} (F^{-1})^T)_{l,m} \int_{\hat{T}} \hat{\varphi}_\alpha \hat{\varphi}_\beta \rmd x
\]
and as $\int_{\hat{T}} \hat{\varphi_\alpha} \hat{\varphi_\beta}\rmd x = \delta_{\alpha,\beta}$ 
the matrix is $(3 \times 3)$-block-diagonal and the inversion is trivial. The
computational effort is obviously of order $\mathcal{O}(N)$ where~$N$ is the
number of unknowns.

\subsubsection{Curved elements}
If we consider curved elements or non-constant material
parameters~$\varepsilon$ and~$\mu$, the approach 
has to be modified as the mass matrix arising from~\eqref{massintegral} may be 
fully occupied. Let's go a step back and consider a similar scalar 
problem\footnote{extensions to 3D are straightforward} with a non-constant
coefficient~$\varepsilon$:
\begin{eqnarray*}
\text{Given:} & & f(v) = \int_T fv \rmd x\\
\text{Find:} & & u, \text{ s.t. } \int_T \varepsilon uv\rmd x = \int_T fv\rmd x 
\end{eqnarray*}
We now transform back to the reference element $\hat{T}$
and get
\[
  \int_T \varepsilon uv\rmd x = \int_{\hat{T}} \vert \det(F) \vert \, \varepsilon uv \rmd x =  \int_{\hat{T}} \hat{u} \tilde{v}\rmd x
\]
where $\tilde{v} = \vert \det(F) \vert \, \varepsilon \hat{v}$. If we 
now approximate~$\tilde{v}$ with the same basis we used for~$v$ before,
the mass matrix is diagonal again. Nevertheless the evaluation of the 
functional~$f(v)$ has to be transformed as well:
\[
  \int_T fv\rmd x = \int_T \frac{1}{\vert \det(F) \vert \, \varepsilon} f\tilde{v} \rmd x = 
  \int_{\hat{T}} \frac{1}{\varepsilon} f \tilde{v} \rmd x
\]
To evaluate the last term we will use numerical integration. But as (in our application)~$f$ is not given pointwise, but in a modal sense, we have to calculate a pointwise 
representation for the numerical integration of $\int_T fv\rmd x$ first:
\begin{eqnarray*}
  \text{Given:} & & f(v) = \int_T fv\rmd x = \int_{\hat{T}} \vert \det(F) \vert \, f\hat{v} \rmd x \\
  \text{Find:} & & f_i, \text{ s.t. } \int_T fv\rmd x = \sum_i \vert \det(F) \vert (x_i) f_i \omega_i v(x_i)
\end{eqnarray*}
Then we can divide (on each integration point) by~$\varepsilon$ and
with those new coefficients we can, by numerical integration, get a
good approximation to $\int_{\hat{T}} \frac{1}{\varepsilon} f
\tilde{v} \rmd x$.  The ``reverse numerical integration'' and the
numerical integration used here can be accelerated by the use of the
sum factorization technique. Doing so the complexity of both ``reverse
numerical integration'' and the numerical integration is
$\mathcal{O}(p^4)$, where $p$ is the polynomial degree. Note that the
approximate inverse $\tilde M_\varepsilon^{-1}$ obtained by this
method is still symmetric and positive definite.

\subsection{Overall computational effort}
In the previous sections we saw that the overall computational effort 
scales linearly with the degrees of freedom~$N$ as long as the elements 
are flat and coefficients are piecewise constant. 
Even for curved elements (and variable coefficients) the 
computational effort is only of order $\mathcal{O}(N^{\frac{4}{3}})$. 
Furthermore no element 
matrices have to be stored. Only the geometric transformations 
and the local topology have to be kept in the memory.

\subsection{Timings}
Let's also state some exemplary numbers that were achieved for this method
and its implementation on an Intel Xeon CPU 5160 at $3.00$ GHz (64 bit) 
(single core) for a tetrahedral mesh with 2078 elements. 
The costs for one step of the symplectic Euler method per 6 scalar degrees 
of freedom are listed in Table \ref{table:timings}.

\begin{table}[ht]
\centering
\subfigure {
\begin{tabular}{c|c}
order $p$ & time $[\mu sec]$ \\
\hline
1 & 0.61 \\
2 & 0.58 \\
3 & 0.71 \\
4 & 0.79 \\
5 & 1.16 \\
6 & 1.24 \\
7 & 1.32 \\
8 & 1.53 \\
9 & 1.66 \\
10 & 1.74 \\
\end{tabular} 
}
\subfigure {
\begin{tabular}{c|c}
order $p$ & time $[\mu sec]$ \\
\hline
1 & 4.89 \\
2 & 2.54 \\
3 & 1.93 \\
4 & 1.79 \\
5 & 2.06 \\
6 & 2.17 \\
7 & 2.33 \\
8 & 2.67 \\
9 & 2.88 \\
10 & 3.04 \\
\end{tabular}
}
\caption{Timings for flat elements (left), using $\mathcal{O}(1)$ floating point operations per dof and 
 curved elements (right) using $\mathcal{O}(p)$) floating point
 operations per dof.}\label{table:timings}
\end{table}

\section{Symbolic derivation of relations}\label{secSymbolics}

In this section we want to describe the symbolic methods that were
employed for finding the desired relations for the polynomial shape
functions. These relations allow for efficient computation of the
discrete curl operations and traces as described in
Section~\ref{secCurls}.  They have been computed by following the
holonomic systems approach~\cite{Zeilberger90,Chyzak00,Koutschan09},
which works for all functions that satisfy sufficiently many linear
differential equations or recurrences or mixed ones;
these relations have to have polynomial coefficients. A large class of
functions (like rational or algebraic functions, exponentials,
logarithms, and some of the trigonometric functions) as well as a
multitude of special functions is covered by this framework. Part of
it are algorithms for the ``basic arithmetic'' (that we will refer to
as ``closure properties''), i.e., given two implicit descriptions for
functions~$f$ and~$g$, respectively, we can compute such descriptions
for $f+g$, $fg$, and for functions obtained by certain substitutions
into $f$ or $g$. All computations in this section have been performed in 
Mathematica using our package \texttt{HolonomicFunctions} (it is freely
available from the website http://www.risc.uni-linz.ac.at/research/combinat/software/).

\subsection{Introductory example}\label{secSymbolicExample}
For demonstration purposes we show how to derive automatically the
rewriting formula~\eqref{rewriteLegendre} for Legendre polynomials~$P_n(x)$.
It is well known that these orthogonal polynomials satisfy some linear 
relations, e.g., the second order differential equation
\[
  (x^2-1)P''_n(x) + 2xP'_n(x) - n(n+1) P_n(x) = 0
\]
or the three term recurrence
\[
  (n+2)P_{n+2}(x) - (2n+3)xP_{n+1}(x) + (n+1)P_n(x) = 0.
\]
We will represent such
linear relations in the convenient operator notation, using the
symbols $\D{x}$ for the partial derivative with respect to~$x$, and
$\S{n}$ for denoting the shift operator with respect to~$n$. Then the two
relations above are written as
\[
  (x^2-1)D_{\!x}^2+2xD_{\!x}+(-n^2-n)
\]
and
\[
  (n+2)S_{\!n}^2+(-2nx-3x)S_{\!n}+(n+1),
\]
respectively, and we identify operators and relations with each other.
The operators can be regarded as elements of a (noncommutative)
polynomial ring in $\S{n}$ and $\D{x}$ with coefficients being
rational functions in $\Q(n,x)$.  We can obtain additional relations
for $P_n(x)$ by combining the given relations linearly, or by shifting
and differentiating them. In the operator setting these operations
correspond to addition and multiplication (from the left) and we can
refer to the set of all operators obtained in this way as the
annihilating left ideal generated by the initially given operators.
In the following we will represent annihilating ideals by means of
their Gr\"obner bases; these are special sets of generators that allow
for deciding the ideal membership problem (i.e., the question whether
some relation is indeed valid for the function under consideration)
and for obtaining unique representatives of the residue classes modulo
the ideal (see~\cite{Buchberger65}).  All algorithms mentioned below
will require Gr\"obner bases as input.  A Gr\"obner basis of the
annihilating ideal of the Legendre polynomials is given by
\[
  G = \big\{(n+1)S_{\!n}+(1-x^2)D_{\!x}+(-nx-x), (x^2-1)D_{\!x}^2+2xD_{\!x}+(-n^2-n)\big\}.
\]

Our main task will be to find elements with certain properties in an
annihilating ideal; this can be done via an ansatz as we demonstrate
now. The relation~\eqref{rewriteLegendre} that we are going to recover
connects $P'_{n+2}(x)$, $P'_n(x)$, and $P_{n+1}(x)$, and its coefficients
are free of~$x$. These facts translate to an ansatz operator of the form
\[
  A = c_1(n)\D{x}\S{n}^2 + c_2(n)\D{x} + c_3(n)\S{n}
\]
where the coefficients~$c_i$ are rational functions in $\Q(n)$, and
hence free of~$x$ as required. We have to determine the $c_i$ such
that the operator~$A$ is an element of the left ideal~$I$ 
generated by~$G$, so that $A(P_n(x))=0$.  For this purpose
we use the Gr\"obner basis~$G$ to compute the unique representation of
the residue class of $A$ modulo~$I$ (it is achieved by reduction). We
have $A\in I$ if and only if the residue class is represented by the
zero operator and hence we can equate all its coefficients to zero,
obtaining the following two equations
\begin{eqnarray*}
  c_1(2nx^2+3x^2-n-2) + c_2(n+1) + c_3(x^2-1) & = & 0,\\
  c_1(n+1)(2n+3)x + c_3(n+1)x & = & 0.
\end{eqnarray*}
Note that in these equations the variable~$x$ occurs, since it is contained in the
coefficients of~$G$. We get a solution that is
free of $x$ by performing a coefficient comparison with respect to this variable.
This yields in the end the linear system 
\[
  \begin{pmatrix}-n-2 & n+1 & -1\\ 2n+3 & 0 & 1\\ (n+1) (2 n+3) & 0 & n+1\end{pmatrix}
  \begin{pmatrix}c_1\\ c_2\\ c_3\end{pmatrix} = 0
\]
whose solution is
\[
  c_1 = -1,\quad c_2 = 1,\quad c_3 = 2n+3,
\]
and this gives rise to the desired relation.

Now what do we do if we don't know the exact shape of the ansatz as
given here by~$A$?  Then we have to include all possible monomials
$\D{x}^i\S{n}^j$ up to some total degree into our ansatz. Looping over
the degree, we will finally find the relation, but the effort can be
tremendous. Therefore, as a preprocessing step, we determine the shape
of the ansatz by modular computations. This means plugging in concrete
values for some of the variables and reducing all integers in the
coefficients modulo some prime. These techniques have been described
in detail in~\cite{Koutschan09} and they are crucial for
getting results in a reasonable time.

All these steps have been implemented in the package
\texttt{HolonomicFunctions} and it computes the relation~\eqref{rewriteLegendre}
immediately:
\begin{mma}
 \In <\!< \ |HolonomicFunctions.m|\\
 \Print \parbox{0.92\textwidth}{
   \vskip 1ex
   \fbox{\parbox{0.9\textwidth}{
     HolonomicFunctions package by Christoph Koutschan, RISC-Linz,\\ 
     Version 1.3 (25.01.2010)\\ 
     $\longrightarrow$ Type ?HolonomicFunctions for help}}
   \vskip 1ex}\\
 \In |FindRelation|\big[|Annihilator|[|LegendreP|[n,x]],\:|Eliminate|\to x\big]\\
 \Out \{S_{\!n}^2D_{\!x}+(-2n-3)S_{\!n}-D_{\!x}\}\\
\end{mma}

\subsection{Relations for the shape functions}

A core functionality of our package
\texttt{HolonomicFunctions}~\cite{Koutschan09} is to execute closure
property algorithms (e.g., for addition, multiplication, and
substitution) on functions represented by their annihilating ideals.
We can now use these algorithms to obtain annihilating ideals for
the shape functions $\varphi$, since their definition in terms of
Jacobi and Legendre polynomials involves just the above mentioned operations.

\subsubsection{The 2D case}\label{secSymb2D}
We first consider triangular finite elements in two dimensions.  For
these, the shape functions are defined as in~\eqref{phi2D}.
Analogously to the one-dimensional example in
Section~\ref{secGradients} we want to express the partial derivatives
(with respect to $x$ and $y$, respectively) in terms of the original
shape functions.  So the goal is to find relations (free of $x$ and~$y$)
that connect the partial derivatives with the original function. More
concretely, we are looking for a relation that allows to express some
linear combination of shifts of
$\frac{\mathrm{d}}{\mathrm{d}x}\varphi_{i,j}(x,y)$ as a linear
combination of shifts of $\varphi_{i,j}(x,y)$ (and similarly
for~$y$). This corresponds to an operator of the form
\begin{equation}\label{phiOp}
  \sum_{(m,n)\in\N^2}c_{1,m,n}(i,j)\D{x}\S{i}^m\S{j}^n+
  \sum_{(m,n)\in\N^2}c_{0,m,n}(i,j)\S{i}^m\S{j}^n
\end{equation}
where the yet unknown coefficients $ c_{d,m,n}\in\Q(i,j)$ do
not depend on $x$ and $y$, and the sums have finite support.

Since we have to find such a relation in the annihilating ideal for
$\varphi_{i,j}(x,y)$, it is natural to start by computing a Gr\"obner basis for this ideal.
The package \texttt{HolonomicFunctions} provides a command \textbf{Annihilator}
that analyzes a given mathematical expression and performs the necessary
closure properties for obtaining its annihilating ideal. So in our example
we can just type
\begin{mma}
 \In |ann|=|Annihilator|[(1-x)\pow i*|LegendreP|[i,2y/(1-x)-1]*\linebreak
   \qquad\qquad |JacobiP|[j,2i+1,0,2x-1],\{|S|[i],\,|S|[j],\,|Der|[x],\,|Der|[y]\}];\\
\end{mma}
\noindent and after a second we have the result (which is
already respectable in size, namely 340kB, corresponding to about 10
pages of output).

Having implemented noncommutative Gr\"obner bases, our first attempt
was to use them for eliminating the variables $x$ and $y$. But it soon
turned out that this attempt did not produce optimal results, and in
addition the computations were very time-consuming. Therefore we came up
with the ansatz described in Section~\ref{secSymbolicExample}. We use it now to
compute the desired relations (both computations take less than a minute):
\begin{mma}
 \In |FindRelation|[|ann|,|Eliminate|\to\{x,y\},|Pattern|\to\{\blank,\blank,0\mid 1,0\}]\linebreak // |Factor|\\
 \Out \{(2i+j+5)(2i+2j+5)S_{\!i}S_{\!j}^2D_{\!x}+(j+3)(2i+2j+5)S_{\!j}^3D_{\!x}+\linebreak
      \phantom{\{}2(2i+3)(i+j+3)S_{\!i}S_{\!j}D_{\!x}-2(2i+1)(i+j+3)S_{\!j}^2D_{\!x}-\linebreak
      \phantom{\{}2(i+j+3)(2i+2j+5)(2i+2j+7)S_{\!i}S_{\!j}-(j+1)(2i+2j+7)S_{\!i}D_{\!x}-\linebreak
      \phantom{\{}2(i+j+3)(2i+2j+5)(2i+2j+7)S_{\!j}^2-(2i+j+3)(2i+2j+7)S_{\!j}D_{\!x}\}\\
 \In |FindRelation|[|ann|,|Eliminate|\to\{x,y\},|Pattern|\to\{\blank,\blank,0,0\mid 1\}]\linebreak // |Factor|\\
 \Out \{(2i+j+6)(2i+j+7)(2i+2j+7)S_{\!i}^2S_{\!j}^2D_{\!y}-(j+3)(j+4)(2i+2j+7)S_{\!j}^4D_{\!y}-\linebreak
      \phantom{\{}4(j+2)(i+j+4)(2i+j+6)S_{\!i}^2S_{\!j}D_{\!y}+4(j+3)(i+j+4)(2i+j+5)S_{\!j}^3D_{\!y}+\linebreak
      \phantom{\{}(j+1)(j+2)(2i+2j+9)S_{\!i}^2D_{\!y}-4(2i+3)(i+j+4)(2i+2j+7)(2i+2j+9)S_{\!i}S_{\!j}^2-\linebreak
      \phantom{\{}(2i+j+4)(2i+j+5)(2i+2j+9)S_{\!j}^2D_{\!y}\}\\
\end{mma}
Here the option \textbf{Pattern} specifies the admissible exponents for
the operators, e.g., in the first case we allow any exponent for the 
shift operators, whereas~$\D{x}$ may occur with power at most~$1$ only, 
and~$\D{y}$ must not appear at all in the result.

\subsubsection{The 3D case}

When dealing with tetrahedra in three dimensions, 
the shape functions are denoted by $\varphi_{i,j,k}(x,y,z)$ and 
are defined by
\[
  %\varphi_{i,j,k}(x,y,z)=
  (1-x-y)^i (1-x)^j 
  P_i\!\left(\textstyle\frac{2z}{1-x-y}-1\right) 
  P_j^{(2i+1,0)}\!\left(\textstyle\frac{2y}{1-x}-1\right) 
  P_k^{(2i+2j+2,0)}(2x-1).
\]
Again they have the nice property of being $L^2$-orthogonal on the reference
tetrahedron
\[
  T=\{(x,y,z)\in\R^3\mid x\geq 0\land y\geq 0\land z\geq 0\land x+y+z\leq 1\}.
\]

Computing an annihilating ideal for $\varphi_{i,j,k}(x,y,z)$ 
is already much more involved than in the 2D case:
\begin{mma}
 \In |phi|=(1-x-y)\pow i\,(1-x)\pow j\,|LegendreP|[i,2z/(1-x-y)-1]\linebreak
     |JacobiP|[j,2i+1,0,2y/(1-x)-1]\,|JacobiP|[k,2i+2j+2,0,2x-1];\\
 \In |Timing|[|ann|=|Annihilator|[|phi|,\:\{|Der|[x],\,|S|[i],\,|S|[j],\,|S|[k]\}];]\\
 \Out \{359.686,|Null|\}\\
\end{mma}
% Timings on gonzales:
% Der[x]: 359.686s, 117MB.  Der[y]: 152.286s, 67MB.  Der[z]: 44.2228s, 16MB.
\noindent The Gr\"obner basis for this annihilating ideal is about
117MB in size (corresponding to several thousand of printed
pages). Note also that it is more efficient to consider only one
derivation operator, and compute annihilating ideals for each of the
cases $\frac{\mathrm{d}}{\mathrm{d}x}$,
$\frac{\mathrm{d}}{\mathrm{d}y}$, and $\frac{\mathrm{d}}{\mathrm{d}z}$
separately (this applies to the 2D case, too).

In principle, the desired relations for the 3D case can be found in the
same way as for two dimensions. As described in Section~\ref{secSymbolicExample}
we find by means of modular computations that the ansatz (for the case
$\frac{\mathrm{d}}{\mathrm{d}x}$) contains the $16$~monomials
\[
  \begin{array}{l}
  S_{\!i}S_{\!j}S_{\!k}^2D_{\!x}, S_{\!i}S_{\!k}^3D_{\!x}, S_{\!j}^2S_{\!k}^2D_{\!x}, S_{\!j}S_{\!k}^3D_{\!x},
  S_{\!i}S_{\!j}S_{\!k}D_{\!x}, S_{\!i}S_{\!k}^2D_{\!x}, S_{\!j}^2S_{\!k}D_{\!x}, S_{\!j}S_{\!k}^2D_{\!x}, \\
  S_{\!i}S_{\!j}S_{\!k}, S_{\!i}S_{\!j}D_{\!x}, S_{\!i}S_{\!k}^2, S_{\!i}S_{\!k}D_{\!x}, S_{\!j}^2S_{\!k}, 
  S_{\!j}^2D_{\!x}, S_{\!j}S_{\!k}^2, S_{\!j}S_{\!k}D_{\!x}.
  \end{array}
\]
However, in order to compute the corresponding coefficients, we did not succeed
with the standard approach used in Section~\ref{secSymb2D}. Instead, we had to
employ modular techniques again for many interpolation points, and then interpolate
and reconstruct the solution.

\section{Conclusion}

We have presented an efficient implementation for solving the
time-domain Maxwell's equations with a finite element method that
uses discontinuous Ga\-ler\-kin elements. Besides many other optimizations
that speed up the whole simulation, the usage of certain recurrence
relations for the shape functions allows for a fast evaluation of
gradients and traces.  These relations have been derived symbolically
with computer algebra methods.

It is widely believed that the mathematical subjects ``numerical
analysis'' and ``symbolic computation'' do not have much in common, or
even that they are kind of orthogonal.  Experts from both areas can
barely communicate with each other unless they don't talk about work.
It was the great merit of the project SFB F013 ``Numerical and
Symbolic Scientific Computing'' that had been established in 1998 at
the Johannes Kepler University of Linz, Austria, to bring together
these two communities to identify potential collaborations.  We
consider our results as a perfect example for such a fruitful
cooperation.

\paragraph{Acknowledgement}
We would like to thank Veronika Pillwein for making contact between
the first- and the last-named author and for kindly supporting our
work by interpreting between the languages of symbolics and numerics.

%\bibliographystyle{plain}
%\bibliography{/home/ckoutsch/tex/literatur.bib}

\end{document}